\documentclass[12pt,leqno,twoside]{amsart}

\let\mgp=\marginpar
\makeatletter\if@twoside
  \def\marginpar#1{}
\else
  \usepackage{showlabels} \def\marginpar#1{\mgp{\raggedright \tiny #1}}
\fi\makeatother

\usepackage{times} 
\usepackage{latexsym}
\usepackage{verbatim}
\usepackage{graphicx,url,enumerate}
\usepackage{amssymb,amsfonts,amsmath}
\usepackage[section]{placeins}

\makeatletter\if@twoside
\else
\fi\makeatother
\newtheorem{theorem}{Theorem}[section]
\newtheorem*{theorem*}{Theorem}
\newtheorem{corollary}[theorem]{Corollary}
\newtheorem{lemma}[theorem]{Lemma}
\newtheorem{proposition}[theorem]{Proposition}

\def\ndthm{Nondegeneracy Theorem}
\newtheorem*{nondegthm}{\ndthm}
\def\dfthm{Diffeomorphism Theorem}
\newtheorem*{diffeothm}{\dfthm}

\theoremstyle{definition}
\newtheorem*{definition}{Definition}
\newtheorem*{remark}{Remark}
\newtheorem*{example}{Example}

\numberwithin{equation}{section}

\def\kgb/{Gro{\ss}e-Brauckmann}

\renewcommand{\d}{\partial}

\renewcommand{\Im}{\operatorname{Im}}

\renewcommand{\epsilon}{\varepsilon}
\newcommand{\eps}{t} 
\newcommand{\ddtz}{\left.\frac{d}{d\eps}\right|_{\eps=0}} 
\renewcommand{\i}{{\mathbf{i}}}  
\renewcommand{\j}{{\mathbf{j}}}  
\renewcommand{\k}{{\mathbf{k}}}  
\renewcommand{\u}{{\mathbf{u}}}  
\renewcommand{\v}{{\mathbf{p}}}  
\newcommand{\vv}{{\mathbf{v}}}  
\newcommand{\w}{{\mathbf{w}}}  
\newcommand{\isom}{\cong}
\newcommand{\ident}{\equiv}
\newcommand{\after}{\circ}
\let\emptyset=\varnothing

\newcommand{\A}{{\mathcal{A}}}
\newcommand{\B}{{\mathcal{B}}}
\newcommand{\D}{{\mathcal{D}}}
\newcommand{\J}{{\mathcal{J}}}

\renewcommand{\L}{{\mathcal{L}}}
\newcommand{\M}{{\mathcal{M}}}
\newcommand{\Mgk}{{\M'_{g,k}}}
\newcommand{\Mzk}{{\M'_k}} 

\newcommand{\W}{{\mathcal{W}}}

\newcommand{\Pik}{\Pi_\k} 

\def\wtilde#1#2#3{{\mskip#1mu\widetilde{\mskip-#1mu #3\mskip-#2mu}\mskip#2mu}}
\def\ptm{\wtilde41M{}^+}
\def\tV{{\wtilde33V}}

\def\tg{{\wtilde33\gamma}}
\def\tf{{\wtilde52f}{}}
\def\tfi{{\tf^{-1}}}
\def\tnu{{\wtilde32\nu}}

\def\wbar#1#2#3{{\mskip#1mu\overline{\mskip-#1mu #3\mskip-#2mu}\mskip#2mu}}
\def\bSigma{{\wbar22\Sigma}}
\def\bW{{\wbar32W}}
\def\btV{{\wbar22\tV}}
\def\ellf#1{{\wbar1{-1}\ell}_{#1}}
\def\barr{{\wbar11r}}
\def\what#1#2#3{{\mskip#1mu\widehat{\mskip-#1mu #3\mskip-#2mu}\mskip#2mu}}
\def\hM{{\what40{\M}}} 
\newcommand{\U}{\mathcal{U}}
\newcommand{\hJ}{\what40{\J}{}}

\renewcommand{\phi}{\varphi}

\newcommand{\R}{{\mathbb{R}}}

\renewcommand{\S}{{\mathbb{S}}}
\newcommand{\Z}{{\mathbb{Z}}}
\renewcommand{\H}{{\mathbb{H}}}

\newcommand{\SO}{\mathsf{SO}}
\newcommand{\so}{\mathsf{so}}
\newcommand{\E}{\mathsf{E}}
\newcommand{\T}{\mathcal{T}}

\newcommand{\twobytwo}[4]{\left(\begin{smallmatrix}#1&#2\\#3&#4\end{smallmatrix}\right)}

\newcommand{\setm}{\smallsetminus}

\begin{document}
\thispagestyle{empty}

\hspace*{-5cm}

\makeatletter\if@twoside\else
\date{\today}
\fi\makeatother

\title{Coplanar $k$-Unduloids are Nondegenerate}
\author[\kgb/]{Karsten \kgb/}
\address{Technische Universit\"at Darmstadt\\ Fachbereich Mathematik (AG 3)\\
Schlossgartenstr.~7\\ 64289 Darmstadt, Germany}
\email{kgb@mathematik.tu-darmstadt.de}
\author[Korevaar]{Nicholas J. Korevaar}
\address{Department of Mathematics\\ University of Utah\\
155 South 1400 East, JWB 233\\ Salt Lake City, UT 84112\\ United States}
\email{korevaar@math.utah.edu}
\author[Kusner]{\\Robert B. Kusner}
\address{Department of Mathematics\\ University of Massachusetts\\
Amherst, MA 01003\\ United States}
\email{kusner@math.umass.edu}
\author[Ratzkin]{Jesse Ratzkin}
\address{School of Mathematical Science, Aras Na Laoi\\ University College Cork\\ 
Cork, Ireland}
\email{j.ratzkin@ucc.ie}
\author[Sullivan]{John M. Sullivan}
\address{Technische Universit\"at Berlin, MA 3--2\\ 
Stra{\ss}e des 17.~Juni 136\\ 10623 Berlin, Germany}
\email{Sullivan@Math.TU-Berlin.DE}


\begin{abstract}
We prove each embedded, constant mean curvature (CMC) surface 
in Euclidean space with genus zero and finitely many
coplanar ends is nondegenerate: there is no nontrivial 
square-integrable solution to the Jacobi equation, the linearization of 
the CMC condition.  This implies that the moduli space of such coplanar 
surfaces is
a real-analytic manifold and that a neighborhood of these in the full
CMC moduli space is itself a manifold.  Nondegeneracy further implies 
(infinitesimal and local) rigidity in the sense that the asymptotes 
map is an analytic immersion on these spaces, and also that the coplanar 
classifying map is an analytic diffeomorphism.

\end{abstract}

\maketitle

\section{Introduction} \label{intro}

The Delaunay unduloids~\cite{D} are surfaces of revolution with 
constant mean curvature $H\ident1$; they are singly-periodic surfaces 
interpolating between a cylinder of diameter~$1$ and a chain of 
unit spheres.
If a finite topology surface~$M$ is properly embedded in~$\R^3$ with 
constant mean curvature (CMC),
Korevaar, Kusner and Solomon~\cite{KKS} proved that
each end of~$M$ is asymptotic to an unduloid,
and if~$M$ is two-ended then it must be an unduloid.
Their results motivate the following global rigidity question:
Do these unduloid asymptotes determine~$M$ uniquely? 
 
To rephrase this question more formally, we define the  \emph{asymptotes map}  
$\A:\hM_{g,k} \rightarrow \U^k$
assigning to any CMC surface~$M$ its $k$ unduloid asymptotes.  Here 
$\hM_{g,k}$ denotes the premoduli 
space of all~CMC surfaces with genus~$g$ and $k$ ends, and $\U:=\hM_{0,2}$ 
is the space of unduloids.  (A precise definition---including a weakened 
notion of embeddedness due to Alexandrov~\cite{A}---is
given in Section~\ref{sec:alex-surfaces}.)
In general, $\hM_{g,k}$ is a real-analytic variety~\cite{KMP}
and $\A$ is a real-analytic map~\cite{Kus-clay}.
Thus a strong form of the above question becomes:
Is the asymptotes map~$\A$ an embedding?

While it appears difficult to answer this question fully,
it follows from \cite{KMP, Kus-clay} that $\A$ is an embedding
in a neighborhood of any CMC surface $M$ which is \emph{nondegenerate}
in the sense that $u\ident0$ is the only $L^2$ solution to the
Jacobi equation
\begin{equation}\label{eq:jacobi}
\L(u) := \Delta u + |A|^2 u = 0.
\end{equation}
Here $\Delta$ is the Laplace-Beltrami operator and
$|A|^2$ 
is the squared length of the second fundamental form of~$M$.
The Jacobi operator~$\L$ is the linearization at~$M$
of the mean curvature operator,
and so (by the implicit function theorem)
near a nondegenerate surface~$M$, the premoduli
space~$\hM_{g,k}$ is a $3k$-dimensional real-analytic manifold~\cite{KMP}.
(It follows that the moduli space $\M_{g,k}$---where surfaces
differing by a rigid motion of~$\R^3$ are identified---has
dimension $3k-6$ near~$M$.)

At a nondegenerate~$M$, we thus get local rigidity, in the sense that
no nearby surfaces have the same asymptotes, and also infinitesimal
rigidity, in the sense that every nonzero Jacobi field on~$M$
has nonzero first-order effect on the asymptotes
(Proposition~\ref{prop:mfld-3k}).

In this paper we consider \emph{coplanar} CMC surfaces, whose
asymptotic axes all lie in a common plane.
Any coplanar CMC surface has mirror symmetry across that plane~\cite{KKS}.
We will call a coplanar CMC surface with genus zero and $k$ ends a
\emph{coplanar $k$-unduloid}, and we let $\Mzk\subset\M_{0,k}$ denote
the moduli space of all coplanar $k$-unduloids.

No classification results are yet available for CMC surfaces with
higher genus or noncoplanar ends, but we do know \cite{GKS1,GKS2}
that $\Mzk$ is homeomorphic to a certain connected
$(2k-3)$--manifold $\D_k$ of spherical metrics on the disk.

Our principal result here is the following:

\begin{nondegthm}
All coplanar $k$-unduloids are nondegenerate.
\end{nondegthm}

Nondegeneracy implies that, in a neighborhood of any coplanar $k$-unduloid,
$\M_{0,k}$ is a real-analytic $(3k-6)$--manifold
and furthermore that $\Mzk$ is an analytic $(2k-3)$--submanifold.
In other words, $\Mzk \subset \M_{0,k}$ has a tubular
neighborhood of dimension $3k-6$, in which $\Mzk$ is
analytically embedded as a submanifold of codimension $k-3$.

Another common application of nondegeneracy is the construction
of new CMC surfaces from old.  Our result gives a wide class of surfaces
to which one can apply gluing constructions \cite{MP,MPP,R}, which 
typically require the summands to be nondegenerate.

There is no example of an Alexandrov-embedded CMC surface
known to be degenerate (other than~$\S^2$).
But without this embeddedness assumption,
there are many examples of degenerate CMC surfaces,
including of course all compact CMC surfaces
as well as unduloids with bubbletons.

Essential to the proof of our theorem is an understanding of the formal 
tangent space at~$M$ to the premoduli space~$\hM_{g,k}$,
which is the space~$\hJ(M)$ of all tempered Jacobi fields
(see Section~\ref{nondegen-background}).
If $M$ is coplanar, the reflection symmetry allows us to
decompose any Jacobi field into even and odd parts,
preserved and reversed (respectively) by the symmetry.
The maximum principle can be used to show~\cite{KKR} that,
up to equivalence, 
the only bounded odd Jacobi field is the one arising
from translation perpendicular to the mirror plane;
this is not in~$L^2$. Thus all $L^2$ fields are even.
We let $\hJ'(M)$ be the space of even, tempered Jacobi fields,
and $\hJ'_0(M)$ be the subspace of $L^2$ fields.  The first 
can be viewed as the subspace tangent at~$M$ to the corresponding 
premoduli space of coplanar CMC surfaces.   

Our proof of the nondegeneracy theorem relies on the 
construction of a formal differential~$\d\Phi$ for the classifying 
map $\Phi\colon\Mzk\to\D_k$ of~\cite{GKS2}.  We prove that $\d\Phi$ is injective
in order to bound the dimensions of~$\hJ'(M)$ and~$\hJ'_0(M)$.
This injectivity also leads to our second main result:

\begin{diffeothm}
The classifying map $\Phi\colon\Mzk\to\D_k$ for coplanar $k$-unduloids
is a real-analytic diffeomorphism. 
\end{diffeothm}

To understand the differential~$\d\Phi$,
let us recall the construction of~$\Phi$:
a coplanar $k$-unduloid $M$ is decomposed by its symmetry plane
into two halves $M^\pm$, each a disk.  The CMC disk~$M^+$ has
a conjugate cousin~$\ptm$, a minimal disk in~$\S^3$ isometric to~$M^+$,
whose boundary lies along $k$ Hopf circles.
Hopf projection immerses this disk into~$\S^2$,
giving the spherical metric $\Phi(M)$ with exactly $k$ completion
boundary points.

The first step in defining $\d\Phi$ is thus a linearization of
the conjugate cousin construction.  This was defined in~\cite{KKR},
and we reinterpret it in Section~\ref{construct-conj}:
Given any Jacobi field on~$M^+$, we get a cousin field on~$\ptm$;
this can be transplanted back to~$M^+$ via the isometry.

Reminiscent of the interchange of Dirichlet and Neumann 
boundary data for conjugate harmonic functions, one might hope
that even fields would conjugate to odd fields and vice versa.
Unforunately, this is not always true, but the conjugate of 
an even field satisfies a slighly weaker condition---we say it is 
\emph{almost odd}. We also define the corresponding notion of
\emph{almost even} fields, and show that a field is almost even if and 
only if its conjugate is almost odd. 

To complete the definition of~$\d\Phi$, as a linear map
on $\hJ'(M)$ modulo Killing fields, we note (Lemma~\ref{lem:hopf-almodd})
that almost odd fields are exactly those which have a well-defined
first-order effect on the Hopf projection of the boundary of~$\ptm$.

One main tool is Theorem~\ref{thm:dPhi-inj}, which says
the map~$\d\Phi$ is injective, implying the upper bound $\dim\hJ'(M)\le2k$.
Key to this proof is a differential version (Lemma~\ref{lem:necksize-chg})
of the fact that the asymptotic necksizes of a coplanar $k$-unduloid $M$ are
readily visible in both $\A(M)$ and $\Phi(M)$.

Finally, the relative index computation of~\cite{KMP}
can be adapted to the equivariant setting;
this shows $\dim\hJ'(M)\ge2k$, with equality if and only
if $M$ is nondegenerate.  Our upper bound shows we do have equality,
completing the proof of the \ndthm.  Nondegeneracy, in turn, implies
that $\Mzk$ is a smooth manifold whose tangent space is the domain
of~$\d\Phi$; thus the injectivity of~$\d\Phi$ also gives the \dfthm.

We expect that our methods might extend to cover the space~$\Mgk$
of coplanar surfaces of higher genus.
New methods, however, not based on the mirror symmetry of
a coplanar surface, would be needed to answer
the question of whether all Alexandrov-embedded
CMC surfaces $M\in\M_{g,k}$ are nondegenerate.

\subsection*{Acknowledgements}
We thank Rafe Mazzeo for many helpful discussions on this project.
We are grateful to the Aspen Center for Physics for hosting two of us
in June 2004 when this work was begun.

\section{Background}\label{sec:bkgrnd}

\subsection{Alexandrov-embedded CMC surfaces of finite topology}
\label{sec:alex-surfaces}

We study constant mean curvature (CMC) surfaces of finite topology.
Here, CMC surfaces are always scaled and oriented (with inward
unit normal $\nu$) to have $H\ident1$,
and a surface of finite topology is
$\Sigma:=\bSigma\setm\{p_1,\ldots,p_k\}$,
where $\bSigma$ is the closed surface of genus~$g$
and each $p_i$ corresponds to an end~$E_i$ of~$\Sigma$.

We are interested in embedded surfaces, but unlike the case of
minimal surfaces, where the maximum principle guarantees that
embeddedness is preserved under deformations, for CMC surfaces
it turns out to be much more natural to consider a slight weakening
of the embeddedness condition.  
A proper CMC immersion $f\colon\Sigma \to \R^3$ of finite topology
is \emph{Alexandrov-embedded} if each end is embedded and
if $f$ extends to a proper immersion $f \colon W \to \R^3$.
Here $W$ is a three-manifold (necessarily a handlebody of genus~$g$)
bounded by~$\Sigma$, and the inward normal~$\nu$ points into~$W$ along~$\Sigma$.

\begin{definition}
Fixing a topological surface~$\Sigma$ of genus~$g$
with $k$ labeled ends $E_1$, \ldots, $E_k$,
we let~$\hM_{g,k}$ denote
the \emph{premoduli space} of complete, Alexandrov-embedded CMC immersions
$f\colon\Sigma\to\R^3$, modulo reparametrizations of the domain
that preserve the labeling.
The quotient of $\hM_{g,k}$ by the action of Euclidean motions
is the \emph{moduli space} $\M_{g,k} := \hM_{g,k} /\E(3)$ of such CMC surfaces 
of genus $g$ with $k$ ends. 
We equip these spaces with the topology
given by Hausdorff distance on compact sets. 
\end{definition}

Delaunay~\cite{D} classified the CMC surfaces of revolution;
the embedded ones, called \emph{unduloids},
are periodic and are determined up to rigid 
motion by their necksize $n \in (0,\pi]$, the length 
of the shortest closed geodesic.  The case $n=\pi$ is the cylinder
(of radius $1/2$, to give $H=1$),
while the limit $n\to0$ approaches a singular chain of unit spheres.

The only Alexandrov-embedded surfaces for $k\le2$
are the sphere and the unduloids \cite{A,M,KKS}.
That is, $\M_{g,k}$ is empty for $g>0$ when $k\le2$,
while $\M_{0,0}\isom\{*\}$, $\M_{0,1}=\emptyset$ and $\M_{0,2}\isom(0,\pi]$.
By contrast, for $k\ge3$,
Kapouleas~\cite{Kap} showed that all $\M_{g,k}$ are nonempty.
Our results here establish that an open set in $\M_{0,k}$ has
the ``expected'' dimension $3k-6$.

We let $\U:=\hM_{0,2}$ denote the premoduli space of unduloids,
and note that an unduloid $U\in\U$ is specified by the following
data: its \emph{axis} (an oriented line in~$\R^3$), its \emph{necksize}
$n\in(0,\pi]$, and (when $n\ne\pi$) its \emph{phase} (the location
of the necks along the axis).  The necksize and phase together
give a point in an open $2$--disk, so $\U$ can be viewed as
a disk-bundle (indeed, a product) over the space of oriented lines in~$\R^3$.
(Compare \cite{Kus-clay}.)

The unduloids are important for us because, as 
we stated in Section \ref{intro}, each end of an 
Alexandrov-embedded CMC surface is asymptotic to some 
unduloid~\cite{KKS}. This motivates us to define the 
asymptotes map $\A$, and to view our Nondegeneracy Theorem as 
an infinitesimal rigidity statement. 

\subsection{Coplanar $k$-unduloids} \label{coplanar-surfaces}
Our arguments apply to surfaces with a certain symmetry.
Alexandrov's reflection technique~\cite{A} can be adapted~\cite{KKS} 
to show that if $M\in\M_{g,k}$ lies in a half-space bounded by a 
plane~$\bar P$,
then it has mirror symmetry across a plane~$P$ parallel 
to~$\bar P$. Furthermore, $P$ cuts~$M$ into two mirror halves~$M^\pm$,
each of which is a graph over a region immersed in~$P$. We conclude 
that the asymptotic axes of such a surface all lie in~$P$, 
and call such a surface \emph{coplanar}. 

We agree to normalize a coplanar surface~$M$
such that $P$ is the horizontal $\i\j$-plane.
We denote by $\hM{}'_{g,k}\subset\hM_{g,k}$ and $\Mgk\subset\M_{g,k}$
the subspaces consisting of coplanar surfaces with this symmetry plane.
The \emph{upper half} of $M\in\Mgk$ is
$M^+:=\{p\in M: \left<p,\k\right> \ge 0\}$.

Because we work with surfaces whose ends are labeled,
the moduli spaces of coplanar surfaces are disconnected.
In this paper, we deal mainly with the case $g=0$,
meaning that~$M$ is a $k$-punctured sphere and~$M^+$ is a closed
disk with $k$ boundary points removed.  (Note that $f$ is a proper
immersion on~$M^+$, whose $k$ ends correspond to those of~$M$.)
The boundary $\d M^+=M\cap P$ consists of $k$ oriented curves~$\gamma_i$,
each connecting a pair of ends.
The space~$\M'_{0,k}$ has $(k-1)!$ components, corresponding to the
different cyclic orderings of the ends.  We focus attention
on one component by assuming the ends of~$M$ are labeled
cyclically, in the sense that $\gamma_i$ runs from~$E_i$ to~$E_{i+1}$.

\begin{definition}
A coplanar CMC surface of genus~$0$, with $k\ge3$ ends
labeled in cyclic order, is called a
\emph{coplanar $k$-unduloid}.
We denote by $\hM_k\subset\hM{}'_{0,k}$ 
and $\Mzk\subset\M'_{0,k}$ the subspaces consisting of coplanar $k$-unduloids.
\end{definition}

The classification result of~\cite{GKS2} shows that
$\Mzk$ is a connected $(2k-3)$--manifold, showing that
$(k-1)!$ is indeed the number of connected components of $\M'_{0,k}$.

When we write a coplanar $k$-unduloid $M\in\Mzk$ as
$M=f(\Sigma)\subset\R^3$, this is a slight abuse of
notation: First, when considering variations, we are often interested
in the particular parametrization~$f$, rather than the equivalence
class under diffeomorphisms of the domain.  Second, although
a nonembedded surface~$M$ cannot be identified with its
image $f(\Sigma)\subset\R^3$, we usually ignore this distinction,
and think of points on~$M$ as being simply points in~$\R^3$.
Finally, we often ignore the distinction between~$M$ and the
abstract surface~$\Sigma$ altogether.

\subsection{Geometry of the three-sphere}\label{sphere-geometry}
Our arguments are based on the construction
of a conjugate cousin, a minimal surface in the three-sphere.
We identify $\S^3$ with the unit quaternions in
$\H=\R[1,\i,\j,\k]\isom\R^4$,
and $\R^3=T_1\S^3$ with the (pure) imaginary quaternions.
As in any Lie group, we can use left-translation to identify
the tangent spaces to~$\S^3$ at any two points.
In particular, a vector $\u\in T_1\S^3\isom\R^3$ is
associated to $p\u\in T_p\S^3$, given by quaternion multiplication.

Important for us will be the \emph{Hopf projection} from~$\S^3$
to~$\S^2$.  Given $\u\in\S^2$, a unit imaginary quaternion,
the $\u$-Hopf projection is defined by $\Pi_\u(p) := p\u p^{-1}$.
The preimage of a point under~$\Pi_\u$ is a great circle in~$\S^3$,
called a \emph{$\u$-Hopf circle}.  The $\u$-Hopf circle through~$p\in\S^3$
has tangent vector~$p\u$ there; in other words, the tangent 
vectors to the $\u$-Hopf circles form the left-invariant vector 
field with value~$\u$ at $1\in\S^3$. 
The collection of all $\u$-Hopf circles,
for fixed~$\u$, foliates~$\S^3$ and is called the $\u$-Hopf fibration.
Usually, we take $\u=\k$; we call $\k$-Hopf circles simply Hopf circles.

We will be interested in the infinitesimal isometries (Killing fields)
on~$\R^3$ and~$\S^3$, because they restrict to Jacobi fields on
CMC surfaces and their conjugates. 
The Killing fields on~$\R^3$ are generated by infinitesimal
translations and rotations.  Given a vector $\u\in\R^3$ the translation
by~$\u$ is the constant vector field $\tau_\u(p)=\u$, while
the rotation around~$\u$ is $\rho_\u(p)=p\u-\u p=-2\u\times p$.
(We include the factor~$2$ because we think of a finite
rotation quaternionically as $p\mapsto e^{-t\u} p e^{t\u}$,
which rotates by angle~$2t$ around the axis~$\u$.)
Similarly, the Killing fields on~$\S^3$ are generated by
infinitesimal left-translations and right-translations.
Given $\u\in\R^3=T_1\S^3$, the left-translation by~$\u$
is the (right-invariant) field $\ell_\u(p)=\u p$,
while the right-translation by~$\u$
is the (left-invariant) field $r_\u(p)=p\u$.  Note that
$r_\u$ is the field of unit tangent vectors to the $\u$-Hopf circles.

\subsection{Conjugate cousins and the classifying map}\label{sec:cousin}
Any simply connected CMC surface in~$\R^3$
has a conjugate cousin minimal surface in~$\S^3$.
We are interested in the case of the upper half~$M^+$
of a coplanar $k$-unduloid $M\in\Mzk$.

We pull back the metric on~$M\subset\R^3$ to~$\Sigma$,
meaning that $f$ is by definition an isometry.
(Indeed, we could identify~$\Sigma$ with~$M$.)
We let $J$ denote a rotation by~$\pi/2$ in the tangent
plane to~$M$ (or equivalently,~$\Sigma$), with sign chosen
so that $\u\times J\u$ is the inward normal~$\nu$ to~$M$.

Then the conjugate cousin of~$M^+$ is a minimal surface
$\ptm=\tf(\Sigma^+)\subset\S^3$.
Here $\Sigma^+:=f^{-1}(M^+)$ and the
immersions~$f$ and~$\tf$ are related by the first-order cousin equation
\begin{equation}\label{eq:cousin}
d\tf=\tf\; df \circ J
\end{equation}
of~\cite{GKS1}.
From this equation, we see that $\tf$ is also an isometry,
and that the normal $\nu$ of $M^+$ left translates to the normal
$\tnu = \tf \nu$ of~$\ptm$. 

As an example, if $M$ is the unit sphere $\S^2=\R^3\cap\S^3$, then $\tf=f$.
(This of course depends on---and indeed characterizes---our sign convention
for~$J$.)

Lawson's original second-order description~\cite{Law}
defined conjugate cousins only up to rigid motion.
But in our first-order approach, given~$f$,
the conjugate~$\tf$ is determined up to a left-translation
in~$\S^3$.  Conversely, given a minimal $\tf\colon \Sigma^+\to\S^3$,
equation~\eqref{eq:cousin} determines~$f$ up to translation in~$\R^3$.

Each of the $k$ boundary curves~$\gamma_i$ of~$M^+$
has constant conormal~$\k$, since it lies in a horizontal plane.
It follows from~\eqref{eq:cousin} that the cousin boundary
curves~$\tg_i$ are Hopf circles. 
Thus under $\k$-Hopf projection,
$\tg_i$ collapses to a point $\v_i\in\S^2$.

On the other hand, the fact that $M^+$ is a graph
means that its normal points downwards: 
$\big<\nu,\k\big> \le 0$, and the strong maximum
principle implies the inequality is strict on the interior.
There we have $\big<\tnu,\tf\k\big>=\big<\nu,\k\big><0$,
which means the interior of the surface $\ptm$ is transverse
to Hopf circles.

It follows that the map $\Pik\after \tf\colon\Sigma^+\to\S^2$
is an immersion on the interior~\cite{GKS2}, so it
induces a spherical metric on the open disk.
The completion boundary of this metric consists exactly
of $k$ points, corresponding to the~$\tg_i$ and developing to the~$\v_i$.
We denote by~$\D_k$ the space of such \emph{$k$-point metrics}:
spherical metrics on the open disk with $k$ completion boundary points.

This discussion has defined a real-analytic map~$\Phi\colon\Mzk\to\D_k$,
taking~$M$ to the metric induced by $\Pik\after\tf$.
The main result of~\cite{GKS2} is that $\Phi$ is a homeomorphism,
explaining why we call it the \emph{classifying map}.

The space~$\D_k$ projects onto a certain space~$\T_k$
of $k$-tuples in~$\S^2$ modulo rotation, by mapping a metric
$D\in\D_k$ to the images~$\v_i$ of its $k$ completion points
(under a developing map to the sphere).
The space~$\T_k$ consists of $k$-tuples where
consecutive points (in cyclic order) are distinct,
and furthermore omits $k$-tuples of the form $(p,q,\ldots,p,q)$
which use only two points of~$\S^2$.  Since rotations act
freely on such $k$-tuples, $\T_k$ is a $(2k-3)$--manifold.
The projection $\D_k\to\T_k$ is a local diffeomorphism~\cite{GKS2};
a $k$-point metric is determined by its completion points
together with some extra combinatorial data. 
For local considerations, it therefore suffices
to consider just the boundary points $\v_i$
instead of the spherical metric $\Phi(M)$.
In particular, when we study the differential of the classifying map,
we will make use of the fact that
\begin{equation}
T_{\Phi(M)}\D_k \isom
\big(T_{\v_1}\S^2 \times \cdots \times T_{\v_k}\S^2\big) \big/ \so_3.
\label{eq:TDk}
\end{equation}
(Here, $\so_3=\{\rho_\u\}$ is the space of infinitesimal rotations
of $\S^2\subset\R^3$.)

Important later will be the fact~\cite{GKS1} that the spherical distance
between successive points~$\v_{i-1}$ and~$\v_i$ is the necksize~$n_i$
of the end~$E_i$.  It is easy to verify this for an unduloid
and its classifying $2$-point metric, by explicitly computing the cousin
surface, a spherical helicoid.  The general case then follows
from the fact that $E_i$ is exponentially asymptotic to an unduloid.

\subsection{Jacobi fields and nondegeneracy}\label{nondegen-background}

To understand the local structure of CMC moduli spaces,
we study the linearized problem.
Suppose $M=f(S)$ is a CMC surface, where we write the
domain as~$S$ to emphasize that we might be talking
about a $k$-unduloid or about just its upper half.
If we vary $f$ in a one-parameter family of immersions
$f^\eps\colon S\to\R^3$ with $f^0=f$,
the \emph{variation field} $\dot f := \ddtz f^\eps$
of this family can be thought of as an $\R^3$-valued vector field on~$M$.
(In general, given a one-parameter family~$f^\eps$,
we will use a dot to denote the first derivative
with respect to the variation parameter~$\eps$, evaluated at $\eps = 0$.)

\begin{definition}
An $\R^3$-valued vector field~$V$ on~$M$ is a
\emph{Jacobi field} if, for every one-parameter family
$f^\eps\colon S\to\R^3$ with variation field $\dot f = V\!$,
the mean curvature $H^\eps$ is constant to first order: $\dot H\equiv0$.  
We call a Jacobi field \emph{integrable} if it is the
first variation of a one-parameter family of CMC surfaces.
\end{definition}

Since tangential fields serve merely to reparametrize the surface,
this notion clearly depends only on
the normal part $V^\perp:=\langle V,\nu\rangle$.
Indeed, $V$ is a Jacobi field if and only if $V^\perp$ is
a \emph{scalar Jacobi field}, satisfying the
Jacobi equation~\eqref{eq:jacobi}. 

\begin{definition}
We call a CMC surface $M$ \emph{nondegenerate} if the only 
$L^2$ Jacobi fields on $M$ are tangential, or in 
other words, if the only solution $u\in L^2$ to the Jacobi
equation~\eqref{eq:jacobi} is $u \equiv 0$.
\end{definition}

The appropriate space of Jacobi fields
to study the local structure of the moduli space $\M_{g,k}$,
in particular to apply the implicit function theorem,
consists of fields more general than $L^2$:
\begin{definition}
We call a Jacobi field \emph{tempered} if it has sub-exponential
growth on each end.
\end{definition}
We note that by the asymptotics result of~\cite{KKS} 
every integrable Jacobi field is tempered.

In the premoduli space $\hM_{g,k}$, surfaces
differing merely by reparametrization are identified;
this means that tangential Jacobi fields have no effect
on the element $M\in\hM_{g,k}$.  This leads us to
define $\hJ(M)$ as the quotient of the space of all
tempered Jacobi fields modulo tangential fields.
Clearly this is isomorphic to the space of normal
(or scalar) tempered Jacobi fields, but often it
is more convenient to work with equivalence classes,
rather than having to take normal parts.

In the moduli space $\M_{g,k}$, surfaces are
further identified if they differ by rigid motion; 
Killing fields have no effect on the element $M\in\M_{g,k}$.
Thus we define $\J(M)$ as the quotient of $\hJ(M)$ modulo Killing fields.
For a Jacobi field~$V$, we write $[V]$ for its coset in~$\J(M)$.

\subsection{Even and odd fields} \label{even/odd}
Suppose~$V$ is a smooth $\R^3$-valued vector field
(for instance, a Jacobi field)
along a coplanar $k$-unduloid~$M$.
Writing $\sigma:(x,y,z) \mapsto (x,y,-z)$
for the reflection across the symmetry $\i\j$-plane,
we can decompose the vector field~$V$
into even and odd parts: $V = V_+ + V_-$, where
$$V_\pm(p) := \tfrac12\big(V(p)\pm\sigma V(\sigma p)\big).$$
Now consider the behavior along the symmetry curves~$\gamma_i$
of these even and odd parts~$V_\pm$
and of their conormal derivatives~$\d_n V_\pm$.
We find that $V_+$ and~$\d_n V_-$ are horizontal (perpendicular to~$\k$),
while $V_-$ and~$\d_n V_+$ are vertical (parallel to~$\k$).

We will often consider vector fields~$V$ defined only on the upper half~$M^+$.
The discussion above motivates the following definitions.

\begin{definition}
Given a vector field~$V$ on~$M^+$,
consider its behavior on one of the boundary curves~$\gamma_i$.
We say that~$V$ is \emph{even} along~$\gamma_i$
if $V$ is horizontal 
and $\d_n V$ is vertical 
there.  Similarly, $V$ is \emph{odd} along~$\gamma_i$
if $V$ is vertical and $\d_n V$ is horizontal.
If one of these conditions holds along all $k$ boundary
curves, we simply say $V$ is even or odd, respectively.
\end{definition}

Note that if $V$ is odd along~$\gamma_i$, then its
normal component $V^\perp=\big<V,\nu\big>$ vanishes
(i.e., has Dirichlet boundary data) along~$\gamma_i$.
Similarly, if $V$ is even along $\gamma_i$ then 
$V^\perp$ has Neumann boundary data in the sense that
$\d_n V^\perp=0$ along $\gamma_i$.
Conversely, if a function $u$ is Dirichlet (or Neumann), then
the normal field $V=u\nu$ is odd (or even, respectively).

Now consider a Jacobi field on~$M^+$.  Typically, it
cannot be extended to a Jacobi field on all of~$M$.
But an even or odd Jacobi field on~$M^+$ does extend,
by even or odd reflection, respectively.

Starting with a coplanar surface~$M$, let $\hJ'(M)$
denote the space of even tempered Jacobi fields (modulo
tangent fields). If we vary~$M$ within~$\hM{}'_{g,k}$,
the symmetry plane is fixed, and the first variation~$V$
is even and tempered.  That is, we can view $V\in\hJ'(M)$.
Indeed, within $\hJ'(M)$, the cone of integrable Jacobi
fields is the tangent cone to~$\hM{}'_{g,k}$.
Similarly, $\J'(M) \subset\J(M)$ denotes the subspace consisting
of cosets of even fields (modulo Killing fields).

\begin{remark}
Results of~\cite{KMP} show that $\M_{g,k}$ is locally
a real-analytic variety; these were adapted in~\cite{GKS2}
to show the same is true of~$\Mgk$.  As varieties, these
spaces have formal tangent spaces at every point~$M$.
It is straightforward to check that these formal tangent
spaces are (isomorphic to) $\J(M)$ and~$\J'(M)$, respectively.
\end{remark}

The tangent cone $T_M \M_{g,k}$ to~$\M_{g,k}$ at~$M$
is the cone in~$\J(M)$ consisting of all (cosets of)
integrable Jacobi fields.  Similarly
$$T_M \Mgk = T_M \M_{g,k}\cap \J'(M)$$
consists of the integrable even Jacobi fields.

\section{Cousin Jacobi fields and the differential of the classifying map}\label{conj-jac-fields}

Our first main technical tool is the conjugation of Jacobi
fields introduced in~\cite{KKR}.  We redevelop this
theory with two changes: we use vector-valued Jacobi fields
and we work directly in~$\S^3$.  Then we employ this to compute the (formal) differential of the classifying map $\Phi\colon\Mzk\to\D_k$.

Throughout this section, we fix a coplanar $k$-unduloid $M\in\Mzk$,
and consider the conjugate~$\ptm$ of its upper half~$M^+$.
Thus $\ptm$ is a minimal surface in~$\S^3$.

\subsection{Conjugation of Jacobi fields} 
\label{construct-conj}
A Jacobi field for the minimal surface $\ptm$
means a variation vector field which preserves the
minimality condition $H\ident0$ to first order.
We think of such a field as a map $W\colon \Sigma^+\to\H$
whose value at $x\in \Sigma^+$
is tangent to~$\S^3$ (though not necessarily to~$\ptm$) at~$\tf(x)$.

Now suppose $V$ is an integrable Jacobi field on~$M^+$,
the initial velocity of a one-parameter family~$f^\eps$ of
CMC immersions.  At each~$\eps$, there is a cousin
minimal immersion $\tf^\eps\colon \Sigma^+\to\S^3$, well-defined up
to left-translation.  This one-parameter family has initial velocity
$\tV:=\ddtz \tf^\eps$,
a Jacobi field on~$\ptm$.
We call $\tV$ a \emph{cousin} of~$V$;
it is well-defined up to an infinitesimal left-translation $\ell_\u$.

Differentiating the cousin equation
$d\tf^\eps = \tf^\eps \,df^\eps \circ J^\eps$
at $\eps=0$ gives the first-order linear system
\begin{equation}\label{eq:Vcousin}
d\tV = \tV \,df\circ J + \tf \,dV\circ J + \tf\, df\circ \dot J
\end{equation}
relating $\tV$ and~$V$.  We will discuss the meaning of
$\dot J := \ddtz J^\eps$ below.

We have derived~\eqref{eq:Vcousin} assuming that~$V$ is integrable.
But in fact, we can start with \emph{any} Jacobi field~$V$
and solve~\eqref{eq:Vcousin} to give a cousin~$\tV$.
This was shown in~\cite{KKR} for the special case when $V$ is normal.
The general case then follows by linearity since any tangential field
is integrable.

Alternatively, we can simply repeat the argument of~\cite{KKR}
for the general case:
Any point $x\in \Sigma^+$ has a neighborhood~$N$
on which the initial surface $f|_N$ is stable for the CMC variational
problem of minimizing area with fixed enclosed volume.  On such
a stable CMC disk, a standard implicit function argument shows
that all Jacobi fields are integrable~\cite{W}.  Thus~\eqref{eq:Vcousin}
holds on the neighborhood~$N$ of the arbitrary point~$x$, meaning
that it holds everywhere.

Multiplying~\eqref{eq:Vcousin} on the left by~$\tfi$
and composing on the right with~$J$, we can solve for~$dV$:
\begin{equation}\label{eq:Vcousinv}
dV = -\tfi\big(d\tV\circ J + \tV\,df\big) + df\circ\dot J\circ J.
\end{equation}
Note that in both equations,~\eqref{eq:Vcousin} and~\eqref{eq:Vcousinv},
we have left the~$\dot J$ term implicit.
It can be computed equally well from~$V$ or from~$\tV$.
To see this, we again use the fact that Jacobi fields are locally integrable
to express~$\dot J$ as the rate of change of conformal structure in
a one-parameter family~$f^\eps$, or equivalently in the isometric
family~$\tf^\eps$.

We have remarked that~$V$ determines~$\tV$ up to a left-translation~$\ell_\u$;
similarly~$\tV$ determines~$V$ up to a translation~$\tau_\u$.  These facts
are also easy to check from the equations we have derived.  Using linearity,
we reduce to the case that the starting field vanishes, in which
case $\dot J=0$.  If $V\ident0$, then~\eqref{eq:Vcousin} becomes
$d\tV=\tV\,df\circ J$, which is solved by
$\tV=\ell_\u=\u\tf$, since it then reduces to~\eqref{eq:cousin}.
In the other direction,
if $\tV\ident0$, then~\eqref{eq:Vcousinv} reduces to~$dV=0$,
giving $V=\tau_\u$.

We summarize the discussion above as follows:
\begin{proposition}\label{pr:isomorphism}
The cousin operation $V\mapsto \tV$ on Jacobi fields,
given by integrating~\eqref{eq:Vcousin}, is an isomorphism from
the space of Jacobi fields on~$M^+$ modulo translations to
the space of Jacobi fields on~$\ptm$ modulo left-translations.
The inverse isomorphism is given by integrating~\eqref{eq:Vcousinv}.
\qed
\end{proposition}

\subsection{Examples}
Trivial examples of Jacobi fields on~$M^+$ and $\ptm$ are the tangential fields
and the restrictions of Killing fields on the respective ambient spaces,
as described earlier.  We can explicitly compute their cousins.

\begin{lemma}\label{lem:tang-cous}
If $V$ is a tangential vector field on~$M^+$,
then $\tV=\tf\,J(V)$ is a cousin.
\end{lemma}
\begin{proof}
Because $V$ is tangential, we can pull it back under~$f$,
to get $V=df(X)$ for some tangent vector field~$X$ on~$\Sigma^+$.
We think of~$X$ as the derivative of a one-parameter
family of diffeomorphisms of~$\Sigma$.  Since the cousin
equation~\eqref{eq:cousin} respects reparametrization
in the sense that $\wtilde31{f\!\!\after\!\phi}=\tf\after\phi$,
it follows that $\tV=d\tf(X)$.  Using~\eqref{eq:cousin}
and the fact that $df$ commutes with~$J$, we get
\begin{equation*}
\tV=\tf\,df(J(X))=\tf\,J(df(X))=\tf\,J(V).\qedhere
\end{equation*}
\end{proof}

Here, we have avoided computing $\dot J$ by not using~\eqref{eq:Vcousin}
directly.  When $V$ is a Killing field, we know in advance that $\dot J=0$.
Of course, if $V=\tau_\u\ident\u$ is a translation field,
then $\tV\ident0$ is a cousin (and indeed the left-translations
are the other cousins).  It is more interesting
to consider the rotational Killing fields.

\begin{lemma}\label{lem:rho-r}
The infinitesimal rotation~$\rho_\u$ around axis~$\u$
has the right-translation $r_\u$ as a cousin.
\end{lemma}

\begin{proof}
We want to show that $V=\rho_\u=f\u-\u f$ and $\tV=r_\u=\tf\u$
satisfy~\eqref{eq:Vcousin} with $\dot J=0$.
But
$$dV=df\,\u-\u \,df, \qquad d\tV=d\tf\,\u=\tf\,(df\after J)\,\u,$$
so indeed $d\tV=\tV\,df\after J+\tf\,dV\after J$.
\end{proof}

\begin{example}
Suppose $M=\S^2$ is a sphere, with $f=\tf$,
and suppose~$V=\rho_\u=-2\u\times f$ is a rotational Killing field.
The symmetry of the sphere means that~$V$ is tangential,
so by Lemma~\ref{lem:tang-cous} one cousin is~$\tV=2\tf\,J(f\times\u)$.
Using the facts that for~$\S^2$ we have
$J(\vv)=\vv\times f$ and $f^2=-1$, we can simplify this as follows:
\begin{align*}
\tV&= 2f\big((f\times \u)\times f) = f\big((f\times \u)f-f(f\times \u)\big)\\
   &= \tfrac12 f\big(2f\u f +2\u)=f\u-\u f 
    = \tf \u-\u \tf = r_\u - \ell_\u.
\end{align*}
Of course, by Lemma~\ref{lem:rho-r}, another cousin of~$V=\rho_\u$
is $r_\u$.  These two cousins are not equal, but their difference
is of course a left-translation field, in this case~$\ell_\u$.
\end{example}

\subsection{Transplanting the cousin}

Given an $\R^3$-valued vector field~$V$ on~$M^+$, which we
think of as a function $V\colon \Sigma^+\to\R^3$,
there is a \emph{transplanted} vector field $\tf V$ on~$\ptm$,
whose value $\tf(x)V(x)$ at~$\tf(x)$ (given by quaternion multiplication)
is tangent to~$\S^3$ there.
That is, we use our identification of different tangent spaces to~$\S^3$
via left-translation to give a natural isomorphism $V\leftrightarrow \tf V$
between the space of vector fields on~$M^+$ and the space of
such $T\S^3$-valued fields on~$\ptm$.
We will most often transplant fields from~$\ptm$ to~$M^+$,
so we introduce the following notation:
Given a field~$W$ on~$\ptm$,
we write $\bW:=\tfi W$ for the transplant back to~$M^+$.

One can easily check that $W$ is a Jacobi field for~$\ptm$ if and only if
its transplant $\bW$ is a Jacobi field for~$M^+$.  (See \cite[Lemma~6]{KKR}.)

Independent of the particular surface~$M^+$,
the transplant to~$M^+$ of a right-translation field~$r_\u$ on~$\ptm$
is the corresponding translational field: $\barr_\u=\tau_\u$.
On the other hand, transplantation of the other Killing fields
gives nontrivial examples of Jacobi fields.  Of particular
interest is the transplant $\ellf\u$,
a Jacobi field on~$M^+$ given by 
$$\ellf\u\big(f(x)\big) = \tfi (x)\,\u\,\tf(x) = \Pi_\u\big(\tfi(x)\big).$$
As discussed in \cite[Appendix~A]{KKR}, one can recover this Hopf 
projection of the conjugate cousin surface from spinning 
spheres, giving an alternative description of the CMC condition. 

We now return to a general transplant $\bW=\tfi W$.
Using $d\tfi=-\tfi\,d\tf\,\tfi$ with~\eqref{eq:cousin}, we get
\begin{equation*}
d\bW = d\big(\tfi W\big)=(d\tfi)W+\tfi\,dW
      = -\tfi\,d\tf\,\tfi\,W + \tfi\,dW
      = -(df\after J)\,\bW + \tfi\,dW,
\end{equation*}
which gives
\begin{equation}\label{eq:dbW}
\tfi\,dW = d\bW + (df\after J)\,\bW.
\end{equation}
The derivative~$d\bW$ can also be interpreted as the covariant derivative
of~$W$ with respect to the flat connection on~$\S^3$ given by
left-translation (while $dW$, on the other hand,
is the derivative in~$\H$).

As an example, if $W=\ell_\u=\u\tf$,
then we get $d\ell_\u=\u\,d\tf=\ell_\u(df\after J)$.
Then by~\eqref{eq:dbW} we have
$$d\ellf\u = \ellf\u\,(df\after J) - (df\after J)\,\ellf\u
           = 2\,\ellf\u\times(df\after J),
$$
or equivalently
\begin{equation}\label{eq:dlu}
2\,df\times \ellf\u - d\ellf\u \after J = 0.
\end{equation}

We note that the conjugation operation discussed extensively
in~\cite{KKR} was the map $V\mapsto\btV{}^\perp$,
sending a Jacobi field~$V$ to the normal part of the transplant
$\btV$ of its cousin~$\tV$.
Let us consider this transplanted cousin~$\btV$.  Using~\eqref{eq:dbW},
we can rewrite~\eqref{eq:Vcousinv} as
\begin{equation}\label{eq:dV}
\begin{split}
dV &= df\,\btV - \btV\,df - d\btV\after J + df\after\dot J\after J \\
   &= 2\,df\times\btV - d\btV\after J + df\after \dot J\after J,
\end{split}
\end{equation}
where, in the second line, we have used the fact that
both $df$ and $\btV$ have values in~$\R^3=\Im\H$ to rewrite
the quaternionic commutator as a vector cross product.
Note that when $V=\tau_\u$ is constant (so $dV=0$ and $\dot J=0$)
then the conjugate is $\btV=\ellf\u$; in this case~\eqref{eq:dV}
reduces to~\eqref{eq:dlu}.

\subsection{Almost even and almost odd fields}
\label{almost-even/odd}

Given a coplanar $k$-unduloid~$M$, we want to examine
vector fields~$V$ on the upper half~$M^+$ in terms of
their behavior on the boundary.  Remember that $\d M^+$
consists of $k$ curves~$\gamma_i$, each lying in
the horizontal $\i\j$-plane of symmetry.
The conjugate cousin surface~$\ptm$ in~$\S^3$
also has $k$ boundary curves~$\tg_i$; these
lie along Hopf circles.

In Section~\ref{even/odd}, we defined $V$ to be
\emph{even} if, on the boundary, $V$ is horizontal and
its normal derivative $\d_n V$ is vertical.
We can interpret this more geometrically by flowing
the surface in the direction~$V$.  We see that
even fields are exactly those for which
the following properties hold to first order:
\begin{enumerate}[(a)]
\item each boundary curve remains in the fixed $\i\j$-plane with $M^+$ meeting
that plane perpendicularly; and
\item the family of conormal curves on~$M^+$, meeting that
boundary perpendicularly, is preserved.
\end{enumerate}

The notion of \emph{odd} fields is easiest
to interpret geometrically when transplanted to~$\ptm$.  Suppose~$W$ is a
vector field on~$\ptm$.  We say $W$ is \emph{odd}
if its transplant~$\bW$ is odd, that is,
if $\bW$ is vertical along the boundary and
its normal derivative $\d_n\bW$ is horizontal.
Equivalently, at a point $p\in\d\ptm$,
this means that $W$ is parallel to the $\k$-Hopf circle
through~$p$ (along which the boundary curve lies)
and that $\d_n W$ is perpendicular to that Hopf circle.
The geometric interpretation is that odd fields are exactly
those for which, when we flow~$\ptm$ in the direction~$W$,
the following hold to first order:
\begin{enumerate}[(a*)]
\item each boundary curve remains along its fixed $\k$-Hopf circle; and
\item the family of conormal curves on~$\ptm$, meeting that boundary
perpendicularly, is preserved.
\end{enumerate}

To appreciate conditions (b) and (b*), think of them
as saying that orthogonality is preserved;
this will later allow us to control the $\dot J$ term
in the cousin equations~\eqref{eq:Vcousin} 
and~\eqref{eq:Vcousinv}.

As mentioned in Section~\ref{sec:bkgrnd}, our goal
is to use cousin Jacobi fields to convert even
fields (which we want to understand) to odd fields.
Odd fields are already better understood thanks
to the following result~\cite[Prop.~24]{KKR}:

\begin{lemma}\label{lem:bdd-odd}
Any bounded odd Jacobi field $V\in\hJ(M)$ is (up
to tangential components) a multiple of the vertical
translation~$\tau_\k$.  In particular, there are
no odd $L^2$ Jacobi fields.
\qed \end{lemma}

Unfortunately, the conjugate of an even field is 
in general not odd, but instead satisfies a slightly weaker condition, 
which we define now, motivated by the geometric interpretations above.

\begin{definition}\label{def:almost-even/odd}
A field~$V$ on~$M^+$ is \emph{almost even}
if along each boundary curve~$\gamma_i$ it differs
from an even field by some translation~$\tau_{\vv_i}$;
here we can take $\vv_i$ to be vertical.
A field~$W$ on~$\ptm$ is \emph{almost odd}
if along each boundary curve~$\tg_i$ it differs
from an odd field by some left-translation~$\ell_{\w_i}$;
here we can take $\w_i$ perpendicular to $\Pik(\tg_i)\in\S^2$.
\end{definition}

Whereas an even field (to first order) keeps each boundary curve~$\gamma_i$
in the $\i\j$-plane, an almost even field translates each~$\gamma_i$ to a
parallel plane, with velocity~$\vv_i$.
Similarly, whereas an odd field preserves (to first order)
each Hopf circle $\tg_i$, an almost odd field will
left-translate it to another Hopf circle.
To quantify this, we compute the derivative of Hopf projection.

\begin{lemma}\label{lem:dPi}
Left-translation Hopf-projects to rotation:
For a unit imaginary $\u\in\S^2$ and any $p\in\S^3$,
we have $\Pik(e^{t\u}p)=e^{t\u}\Pik(p)e^{-t\u}$.
Infinitesimally,
$$d_p\Pik\big(\ell_\u(p)\big) = d_p\Pik(\u p)
= 2\u\times\Pik(p)= -\rho_\u(\Pik p).$$
That is, $d\Pik(\ell_\u)=-\rho_\u$, in the sense that
the vector field $\ell_\u$ is $\Pik$-related to $-\rho_\u$.
\end{lemma}
\begin{proof}
The definition $\Pik(p):=pkp^{-1}$ immediately yields
the equation for $\Pik(e^{t\u}p)$.  Differentiating
this at $t=0$ gives
\begin{equation*}
d_p\Pik(\u p)= \tfrac{d}{dt}\bigr|_{t=0}e^{t\u}\Pik(p)e^{-t\u}
        = \u \Pik(p)-\Pik(p)\u = 2\u\times \Pik(p).\qedhere
\end{equation*}
\end{proof}
We note that, according to this formula, the only two $\k$-Hopf circles
remaining fixed under $\ell_\u$ are the ones projecting to~$\pm\u\in\S^2$.
This makes sense, since along these circles,
$\ell_\u=\u p=p\k =r_\k$ is tangent.

\begin{lemma}\label{lem:hopf-almodd}
The almost odd fields on~$\ptm$
are exactly those which have a well-defined first-order effect on
$$\Pik(\d\ptm)= (\v_1,\ldots,\v_k)\in(\S^2)^k;$$
this is given by $\rho_{\w_i}$ at~$\v_i$.
\end{lemma}
\begin{proof}
An odd field preserves the boundary Hopf circles by property (a*) above.
An almost odd field~$W$ differs from an odd field
by $\ell_{\w_i}$ along $\tg_i$.
Thus by Lemma~\ref{lem:dPi} its action on $\v_i:=\Pik(\tg_i)$ is
the rotation $\rho_{\w_i}$.
(Note that the action of $\rho_{\w_i}$ on~$\v_i$ is well-defined, since
different choices of~$\w_i$ differ by multiples of~$\v_i$.)
\end{proof}

For local computations along the boundary curve~$\gamma_i$,
we pick an orthonormal coordinate frame $(\d_t,\d_n)$, where $\d_n=J\d_t$ 
and $\tau:=\d_t f=-\tfi\d_n\tf$ is the horizontal tangent vector
along~$\gamma_i$, 
while $\k=\d_n f=\tfi\d_t\tf$ is the constant vertical conormal.  
In these coordinates, we can
give an alternate characterization of almost even and almost odd fields.

Along a boundary curve, an almost even field differs from
an even one by a constant.  Thus it is clear that
a field~$V$ on~$M^+$ is almost even if and only if,
along the boundary~$\d M^+$, the tangential derivative
$\d_t V$ is horizontal, while the normal derivative $\d_n V$
is vertical.  We think of these conditions as properties of the
vector-valued $1$-form $dV$, and have thus proved the first
half of the following lemma:

\begin{lemma}\label{lem:almostodd}
A vector field~$V$ on~$M^+$ is almost even if and only if,
along the boundary~$\d M^+$, the one form~$dV$ is
horizontal when applied to~$\d_t$ and vertical when applied to~$\d_n$.
Similarly, a vector field~$W$ on~$\ptm$ is almost odd if and only if
the form~$2\,df\times\bW-d\bW\after J$ is
horizontal when applied to~$\d_t$ and vertical when applied to~$\d_n$.
\end{lemma}

\begin{proof}
To prove the second statement, let $\omega:=2\,df\times\bW-d\bW\after J$.
By~\eqref{eq:dlu}, this form~$\omega$ is unchanged
if we add a left-translation~$\ell_\u$ to~$W$.

Therefore, given an almost odd~$W$, we are free to assume it is odd
along any given boundary curve, that is, that $\bW$ is vertical
and $\d_n\bW$ is horizontal.  Then $\omega(\d_t)=2\tau\times\bW-\d_n\bW$
is indeed horizontal, and $\omega(\d_n)=2\k\times\bW+\d_t\bW=\d_t\bW$
is indeed vertical.

Conversely, given~$W$ for which $\omega$ has
the given properties, we are free to assume that $\bW$ is vertical
at some initial boundary point.
At any point where $\bW$ is vertical, $\k\times\bW$ vanishes,
so $\d_t\bW=\omega(\d_n)$ is also vertical.
Thus, in the unique solution to the ODE, we find that $\bW$ stays vertical
along the whole boundary curve.
\end{proof}

\begin{proposition}\label{prop:ale-alo}
Suppose $V$ is a Jacobi field on~$M^+$
and $\tV$ is a cousin field on~$\ptm$.
Then $V$ is almost even if and only if $\tV$ is almost odd.
\end{proposition}
\begin{proof}
Rewriting equation~\eqref{eq:dV} as
$$ dV- \big(2\,df\times\btV - d\btV\after J\bigr) = df\after \dot J\after J,$$
we see by Lemma~\ref{lem:almostodd} that it suffices to prove
that 
$df\after\dot J\after J$ is
horizontal when applied to~$\d_t$ and vertical when applied to~$\d_n$.

When flowing by an even field~$V$ or an odd field~$\tV$,
it follows from properties (a,b) or (a*,b*) above that,
to first order, the frame $(\d_t,\d_n)$ remains orthogonal along~$\gamma_i$.
The same is true if we flow by an almost even or almost odd field,
since locally the only difference is a Killing field.

Thus in this frame, to first order, $J$ flows from
$\twobytwo0{-1}10$ to a matrix of the form $\twobytwo0{-1/a}a0$.
That is, $\dot J=\twobytwo0{\dot a}{\dot a}0$,
implying that $\dot J\after J$ is diagonal.
Thus $df\after \dot J\after J(\d_t)$ is horizontal,
a multiple of~$\tau$, while $df\after\dot J\after J(\d_n)$ is vertical,
a multiple of~$\k$.
\end{proof}

\subsection{The differential of the classifying map}
We now have all the ingredients we need to compute the
formal differential of the classifying map $\Phi\colon \Mzk\to\D_k$.
Recall that $\Phi$ takes~$M$ to the Hopf projection
of the conjugate cousin $\ptm$ of its upper half.
We write $(\v_1,\ldots,\v_k)\in\T_k$ for the boundary
of~$\Phi(M)$.

Given an almost even field~$V$ on~$M^+$,
by Proposition~\ref{prop:ale-alo} any conjugate~$\tV$ is almost odd,
differing from an odd field by some $\ell_{\w_i}$ on each~$\tg_i$.
Thus, by Lemma~\ref{lem:hopf-almodd}, $\tV$ has the well-defined
effect $-\rho_{\w_i}$ on each $\v_i$.  Of course, $\tV$ itself is only
well-defined up to some left-translation~$\ell_u$, but this just
gives a global rotation $-\rho_\u\in\so_3$ of the whole $k$-tuple,
so we end up with a well-defined element of~$T_{\Phi(M)}\D_k$,
using the characterization~\eqref{eq:TDk} of that space.

\begin{proposition}\label{prop:def-dPhi}
The formal differential $\d\Phi:\J'(M) \to T_{\Phi(M)}\D_k$, 
given by 
$$\d\Phi([V]) := ( -\rho_{\w_1}, \dots, -\rho_{\w_k}),$$
is well-defined. Moreover, when $V$ is integrable
(i.e., when $[V]\in T_M\Mzk$) this is---as suggested by the notation---the
derivative of the classifying map~$\Phi$.
\end{proposition}
\begin{proof}
Suppose $[V]\in\J'(M)$ is given, with $V$ an even Jacobi field.
The construction of~$\d\Phi$ parallels that of~$\Phi$, making
the last statement clear by the chain rule.

The only thing left to check is that $\d\Phi([V])$
is independent of the even representative $V\in[V]$.
But this is straightforward: An even Killing field
is either a horizontal translation (with no effect on~$\tV$)
or the rotation~$\rho_\k$; by Lemma~\ref{lem:rho-r},
the latter has as a cousin~$r_\k$,
which is tangent to $\k$-Hopf circles and thus has no effect on the~$\v_i$.
Similarly, an even tangential field is horizontal on~$\d M^+$,
so by Lemma~\ref{lem:tang-cous} it has a tangential cousin
which is vertical (that is, again tangent to $\k$-Hopf circles)
on the boundary, so again has no effect.
\end{proof}

We can easily give a preliminary characterization of the kernel of this map;
this lemma will later be used in showing that in fact the kernel vanishes.
\begin{lemma}\label{lem:ker-dPhi}
For $[V]\in\J'(M)$, we have $\d\Phi([V])=0$ if and only
if $V$ has an odd cousin~$\tV$.
\end{lemma}
\begin{proof}
Any cousin~$\tV$ is almost odd, meaning that it differs from an odd field
by some $\ell_{\w_i}$ along~$\tg_i$.
Thus~$\tV$ is odd if and only if the vectors~$\w_i$ can be chosen
to vanish.  Of course, the definition of $T\D_k$ means
that $\d\Phi([V])=0$ even if the $k$-tuple is rotating,
that is, if all the~$\w_i$ can be chosen equal to some fixed~$\u$.
But then $\tV-\ell_\u$ is another cousin for which $\w_i=0$;
it is thus odd.
\end{proof}

\section{Moduli space theory and the differential of the asymptotes map}
In this section, we review some basic structure results about CMC moduli spaces.
(See~\cite{KMP} and the analogous
discussion of constant scalar curvature metrics in~\cite{MPU}.)
Then we use the linear decomposition lemma to compute the (formal)
differential of the asymptotes map $\A:\hM_{g,k} \rightarrow \U^k$.

\subsection{Delaunay unduloids}
Recall that an unduloid has a conformal parametrization of the form
$$U:(t,\theta)\mapsto x(t)\,\i + r(t)e^{\theta\i}\,\j.$$
Here $U$ is an unduloid with necksize $n$, positioned so that its
axis is the $\i$-axis and one neck lies in the $\j\k$-plane.
Thus $x(0) = 0$ and $r$ assumes its minimum of $n/2\pi$ at $t=0$. 

We let $\eta$ be the Jacobi field arising from differentiation
with respect to the necksize~$n$; we call~$\eta$ the necksize-change field.
(On the cylinder, where $n=\pi$, the space of necksize-change fields
is two-dimensional, since necks can be inserted at any phase along the axis.)

The necksize-change and Killing fields together form
a six-dimensional subspace of $\hJ(U)$.
(Note that $\rho_\i$ is tangential, so does not contribute here.)
The translation fields are bounded (in fact, periodic in~$t$),
while $\eta$ and the rotation fields grow linearly~$t$.
(On the cylinder the translation $\tau_\i$ is also tangential, but
the extra necksize change compensates in the dimension count.
Here, the necksize changes are bounded.)  

Using the rotational symmetry of~$U$ to separate variables
and expand any Jacobi field in a Fourier series, we see that 
the modes of order $0$ and $\pm 1$ correspond to the geometric
motions described above, while the higher Fourier modes 
grow exponentially on at least one end ($t\to \infty$ or $t\to -\infty$).
On the other hand, if $u \in L^2\bigl([0,\infty) \times \S^1\bigr)$
is a Jacobi field, then $u$ must decay exponentially, and 
there is a lower bound (independent of the necksize $n$)
for the rate of exponential decay. 

We can summarize this discussion as follows.
(See also \cite{KKS}, \cite{KMP}, and \cite[Prop.~20]{MP}.)  
\begin{lemma} \label{lem:del-jacobi-fields}
The unduloid~$U$ with necksize~$n$ is nondegenerate.
Indeed, the space~$\hJ(U)$ of tempered Jacobi fields
is six-dimensional, spanned by the necksize-change and Killing fields.
These fields all grow at most linearly in~$t$, and
are all integrable, meaning~$\hJ(U)$ is the tangent
cone $T_{U}\U$ to the premoduli space $\U=\hM_{0,2}$.
The even part $\hJ'(U)=T_U\U'$ is four-dimensional,
spanned by the necksize-change,
the horizontal translations ($\tau_\i, \tau_\j\,$), and the rotation~$\rho_\k$.
\qed\end{lemma}

\subsection{The Jacobi operator on finite-topology CMC surfaces}
Fixing a CMC surface $M\in\hM_{g,k}$, we want to understand
the asymptotic behavior of tempered Jacobi fields on~$M$.
Following \cite[p.~126]{KMP},
choose a representative~$E_i$ for each end that is a normal
graph over one end of an unduloid~$U_i$;
this allows us to identify~$E_i$ with a subset of~$U_i$.
Cover~$M$ with the~$E_i$ plus a compact set~$E_0$ and
let~$\{\phi_i: 0\le i\le k\}$ be a partition of unity subordinate
to this cover.  We define the \emph{deficiency space}
$$\W=\W(M):=\Big\{W=\sum_{i=1}^k \phi_i W_i : W_i\in\hJ(U_i)\Big\},$$
and for a coplanar $M\in\hM{}'_{g,k}$, we further define the \emph{even deficiency subspace}
$$\W':=\Big\{W=\sum_{i=1}^k \phi_i W_i : W_i\in\hJ'(U_i)\Big\},$$
assuming we have chosen the partition of unity symmetrically.

The elements $W\in\W$ are not Jacobi fields, but their normal parts
$\langle W,\nu\rangle$ have the property that,
along each end, $\L\bigl(\langle W, \nu\rangle\bigr)$ decays exponentially
in the coordinate~$t$,
at a rate depending only on the asymptotic necksize of that end.
One may also regard $\W$ as a quotient space of functions on~$M$ which
are asymptotic to some tempered Jacobi field on each end,
modulo functions asymptotic to zero; we use this interpretation
implicitly in the Remark below Propostion \ref{prop:dim-3k2k}.

To distinguish tempered Jacobi fields from
those which grow exponentially on at least one end,
we use weighted Sobolev spaces~$H^s_\delta(M)$.
We say $u\in H^s_\delta$ if $u \in H^s_\text{loc}$ and,
when restricted to any end, $e^{-\delta t} u \in H^s$.
We write the Jacobi operator~$\L$ as
$$\L_\delta\colon  H^{s+2}_\delta \to H^s_\delta$$
to emphasize the choice of weight~$\delta$ for the domain and codomain.

We recall the linear decomposition lemma \cite[Lemma~2.9]{KMP}:

\begin{lemma} \label{lem:lin-decomp}
There exists $\gamma>0$, depending only on the necksizes of~$M$,
such that for any $\delta\in(0,\gamma)$ any function
$u \in H^{s+2}_\delta(M)$ with $\L(u) \in H^s_{-\delta}(M)$
can be decomposed as
$u = \langle W, \nu \rangle + v$,
where $W \in \W$ and $v \in H^{s+2}_{-\delta}(M)$;
in particular, $u$ has at most linear growth along each end. 
\qed \end{lemma}

\subsection{The differential of the asymptotes map}
An immediate corollary 
of the linear decomposition lemma
is the existence of a formal differential~$\d\A$
for the asymptotes map:

\begin{corollary}\label{cor:dA}
The asymptotes map $\A\colon \hM_{g,k} \to \U^k$
has a formal differential~$\d\A$ which is
defined for all tempered Jacobi fields $V\in\hJ(M)$, and which 
coincides with the derivative of $\A$ on integrable fields.
\end{corollary}
\begin{proof}
Given~$V\in\hJ(M)$, Lemma~\ref{lem:lin-decomp} gives $V=\sum\phi_i W_i + v$,
with $v$ exponentially decaying.
This means $V$ is asymptotic to~$W_i\in\hJ(U_i)$ on end $E_i$.
We now set $\d\A(V):=(W_1,\ldots,W_k)$.
When $V$ is integrable, the $W_i$ must agree with the actual
first-order changes of the asymptotes, so this $\d\A$ is indeed
the differential of~$\A$.
\end{proof}

Given a Jacobi field $V\in\hJ(M)$, we now have two ways of
measuring the rate of change it induces in the asymptotic necksize~$n_i$:
we can compute $\d\Phi([V])$ and find the rate of change of
distance from~$\v_{i-1}$ to~$\v_i$ (Proposition \ref{prop:def-dPhi}),
or we can compute~$\d\A(V)$ and look at
the necksize-change component of~$W_i$ (Corollary \ref{cor:dA}).
Clearly the computations agree for integrable fields, since they both
measure the actual change of necksize along any one-parameter
family of CMC surfaces.  In fact, they always agree:

\begin{lemma}\label{lem:necksize-chg}
The rate of change of necksize~$n_i$ under $V\in\hJ(M)$ can be 
computed equivalently from~$\d\Phi$ or from~$\d\A$.
\end{lemma}
\begin{proof}
Since $V$ and $W_i$ are exponentially asymptotic along~$E_i$,
they have the same effect on~$\v_{i-1}$ and $\v_i$.
But $W_i$ is an integrable field on an unduloid~$U_i$,
so its necksize-change component agrees with the rate
of change of distance from~$\v_{i-1}$ to~$\v_i$.
\end{proof}

\subsection{Dimension counting with (coplanar) symmetry}
The linear decomposition lemma applies in particular to tempered Jacobi fields,
showing their normal parts have at most linear growth.
Since no nonzero element of~$\W$ has normal part in~$L^2$, we also see
that $L^2$ normal Jacobi fields must decay exponentially on each end.
More precisely, if we write $\hJ_0(M)$ for the space of $L^2$ Jacobi fields
and $\hJ'_0(M)$ for the subspace of even fields,
then for sufficiently small $\delta>0$, we find that
$\hJ(M)\isom\ker(\L_\delta)$
and similarly $\hJ_0(M)\isom\ker(\L_{-\delta})$.

The duality between~$\L_\delta$ and~$\L_{-\delta}$
suggests that we use the relative index theorem of Melrose \cite[Sec.~6.2]{Mel}
to compute the difference in dimension of~$\hJ(M)$ and~$\hJ_0(M)$.
Indeed, if $M\in\M_{g,k}$ has any finite isometry group $G<\SO_3$
(acting perhaps to permute the ends), we can consider the spaces of
$G$-invariant Jacobi fields, and we get:
\begin{proposition} \label{prop:dim-3k2k}
For any $G$-symmetric $M\in\M_{g,k}$ we have
\begin{equation*}
\dim\hJ^G(M)-\dim\hJ^G_0(M)=\tfrac12\dim\W^G. 
\end{equation*}
\end{proposition}
Here, $\dim\W^G$ can be computed as a sum over a set of representative ends,
inequivalent under~$G$.  The contribution of each end is either~$6$
(for an end in general position), or~$4$ (for an end in a mirror plane)
or~$2$ (for an end along a rotation axis of~$G$).  The first two cases
are already familiar to us from Lemma~\ref{lem:del-jacobi-fields};
in the last case, of course, the only
symmetric perturbations of the end are necksize change and translation 
along the axis.

The case of no symmetry, $G=1$, was proved as \cite[Thm.~2.11]{KMP}.
The general case follows from the same relative index calculation,
simply by restricting throughout to symmetric functions.
(Specific cases have been used, for instance, in \cite[Sect.~4.3]{JP}
and \cite[Props.~2.1,~2.2]{HP}.)
For the case of coplanar surfaces, $G\isom\Z_2$ is mirror symmetry,
and the result becomes
\begin{equation}
\dim\hJ'(M)-\dim\hJ'_0(M)=\tfrac12\dim\W'=2k.\label{eq:dim-2k}
\end{equation}

\begin{remark}
A more conceptual proof of Proposition~\ref{prop:dim-3k2k} begins
by observing that $\W$ is a symplectic vector space with respect to
the Gauss--Green form~\cite{KMP}, and that the bounded nullspace
$\B=\B(M) := \hJ(M)/\hJ_0(M)$ may be regarded as an isotropic
subspace of~$\W$.  The Melrose relative index theorem implies
$\dim\B=3k=\tfrac12\dim\W$ (whether or not $M$ is nondegenerate),
so $\B$ is Lagrangian in~$\W$.  If $G$ acts on~$M$ by isometries,
it induces a symplectic $G$-action on~$\W$ preserving~$\B$.  But it
follows from a symplectic linear algebra lemma \cite[Lemma~6.2]{GKS2}
that the fixed-point set $\W^G$ is a symplectic subspace of~$\W$, and
also that  $\hJ^G(M)/\hJ^G_0(M)=\B^G=\B\cap\W^G$ is Lagrangian in~$\W^G$,
yielding the result.
\end{remark}

The following result (\cite[Thm.~3.1]{KMP}, as reinterpreted in~\cite{Kus-clay})
explains the connection between nondegeneracy of~$M$ and the
regularity of the asymptotes map near~$M$.

\begin{proposition}\label{prop:mfld-3k}
If $M$ is nondegenerate then there exists a neighborhood of~$M$ in~$\hM_{g,k}$ 
which is a manifold of dimension~$3k$,
and on which the asymptotes map
$\A:\hM_{g,k} \rightarrow \U^k$ is an embedding.
\qed \end{proposition}

The analagous statement
for coplanar surfaces appears as~\cite[Thm.~5.2]{GKS2}.
Combining it with equation~\eqref{eq:dim-2k}, we get the following:

\begin{theorem}\label{thm:dim-2k}
For a coplanar $M \in \Mgk$, we have $\dim\hJ'(M)\ge 2k$, 
with equality if and only if $M$ is nondegenerate. Also, if 
$M$ is nondegenerate then, in a neighborhood
of~$M$, the premoduli space~$\hM{}'_{g,k}$ is a real-analytic $2k$-manifold
with tangent space~$\hJ'(M)$, and the moduli space~$\Mgk$
is a $(2k-3)$--manifold with tangent space~$\J'(M)$.
\qed
\end{theorem}

\section{Proofs of the main theorems}\label{main-thm-sec}

The key to our main results is the following theorem,
guaranteeing the injectivity of~$\d\Phi$.

\begin{theorem} \label{thm:dPhi-inj}
The formal differential
 $\d\Phi\colon \J'(M)\to T_{\Phi(M)} \D_k$
of the classifying map~$\Phi$,
as given in Proposition~\ref{prop:def-dPhi}, is injective.
Thus $\dim\J'(M)\le2k-3$, so $\dim\hJ'(M)\le2k$.
\end{theorem}

\begin{proof} 
Suppose $V\in\hJ'(M)$ is an even field with $\d\Phi([V])=0$.
By Lemma~\ref{lem:ker-dPhi}, it has an odd cousin~$\tV$.
We claim the normal part of~$\tV$ is bounded.
Then by Lemma~\ref{lem:bdd-odd}
its transplant is a vertical translation~$\tau_{a\k}$.  Equivalently,
$\tV=r_{a\k}$. (All our computations are up to tangential components.)
But we know that $r_{a\k}$ is a cousin of~$\rho_{a\k}$,
so it follows that $V=\rho_{a\k}+\tau_\u$ for some $\u\in\R^3$.
But this is a Killing field, meaning that $[V]=0\in\J'(M)$.
The dimension bounds follow immediately from the known dimension of~$\D_k$.

To prove the claim, we look at $\d\A(V)=(W_1,\ldots,W_k)$.
Here $W_i\in\hJ'(U_i)$ where $U_i$ is the asymptote of end~$E_i$.
Since $\d\Phi([V])=0$, by Lemma~\ref{lem:necksize-chg},
the fields~$W_i$ include no necksize changes.  Thus each $W_i$
is an even Killing field, in the span of $\{\tau_\i,\tau_\j,\rho_\k\}$.
But we know cousins for the Killing fields: translations have
vanishing cousins, while $\rho_\k$ has cousin~$r_\k$ by Lemma~\ref{lem:rho-r}.
Thus $W_i$ has a bounded cousin~$r_{a\k}$, so $V$ has a cousin
which along~$E_i$ is exponentially asymptotic to this
right-translation, and in particular is bounded.
Of course, since left-translations are bounded, the boundedness of
some cousin along each end implies the boundedness of any cousin
along all ends; in particular, the odd cousin $\tV$ is bounded.
\end{proof}

\begin{nondegthm}
All coplanar $k$-unduloids are nondegenerate.
\end{nondegthm}
\begin{proof}
Combining Theorems~\ref{thm:dim-2k} and~\ref{thm:dPhi-inj}, 
we get $\dim\hJ'(M) = 2k$; the equality in Theorem~\ref{thm:dim-2k}
then implies nondegeneracy.
\end{proof}

\begin{diffeothm}
The classifying map $\Phi\colon\Mzk\to\D_k$ of~\cite{GKS2}
is a real-analytic diffeomorphism. 
\end{diffeothm}
\begin{proof}
Using the \ndthm\ and Theorem~\ref{thm:dim-2k},
we see that~$\Mzk$ is a real-analytic manifold of dimension $2k-3$,
with tangent space~$\J'(M)$.  The classifying map~$\Phi$ is
thus a real-analytic map between manifolds of the same dimension,
and Theorem~\ref{thm:dPhi-inj} says its differential is injective.
The theorem then follows by the real-analytic inverse function theorem. 
\end{proof}

\end{document}